\documentclass[12pt]{article}
\usepackage{amsmath, amssymb, graphics}
\usepackage{amsthm}
\usepackage{graphicx}
\usepackage{color}
\usepackage{authblk}

\newtheorem{thm}{Theorem}[section]
\newtheorem{corollary}{Corollary}[section]
\newtheorem{proposition}{Proposition}[section]
\newtheorem{lemma}{Lemma}[section]
\newtheorem{example}{Example}[section]
\newtheorem{defn}{Definition}[section]
\newtheorem{remarks}{Remarks}[section]

\def\N{\mathbb{N}}

\def\R{\mathbb{R}}

\def\Var{{\rm Var} \,}
\def\ee{\varepsilon}
\def\E{{\mathbb E}}
\def\cE{{\mathcal E}}
\def\P{{\mathbb P}}
\def\Cox{\hfill \Box}

\def\disp{\displaystyle}
\def\one{{\bf 1}}
\def\|{{\, | \, }}
\def\F{{\mathcal F}}
\def\G{{\mathcal G}}
\def\B{{\mathcal B}}

\def\Y{{\bf Y}}

\def\m{{\bf m}}
\def\tree{{\mathcal T}}

\def\pace{{\mathcal J}}
\def\zerotree{{\emptyset_*}}
\def\diam{{\rm diam}\,}
\def\nt{{\mathfrak n}}
\def\vertex{\Xi}
\def\rt{{\bf 0}}
\def\parent{{\rm par}\,}
\def\ulam{{\mathcal U}}
\def\VU{{\bf V}}  
\def\bfa{{\bf a}}
\def\timesf{{\mathcal Y}}

\title{Euclidean Embedding of the \\Poisson Weighted Infinite Tree\\ and Application to Mobility Models} 

\author[1]{R. W. R. Darling } 
\affil[1]{National Security Agency\footnote{Mathematics Research Group,  Fort George G. Meade, MD 20755-6484}}

\author[2]{Robin Pemantle 
} 
\affil[2]{University of Pennsylvania\footnote{Mathematics Dept., 209 South 33rd Street, Philadelphia, PA 19104-6395}}

\date{\today}

\begin{document}
\maketitle
\begin{abstract}
Continuous time branching models are used to create random fractals in a Euclidean space, whose Hausdorff dimension is controlled by an input parameter. Finite realizations are applied in modelling the set of sites visited in models of human and animal mobility.\\


\end{abstract}

\section{Introduction and motivation}

\subsection{Context: stochastic mobility models}
Zoologists, social scientists, and online advertisers which use locations inferred from communication metadata to target content, are building algorithms which process data sets containing sporadic reports of animal or human mobility. Such data sets are expensive to create by scientific means, and their distribution is limited by privacy and security concerns. This leads to the need to create realistic synthetic data for algorithm testing.

After examining data sets from several sources, the authors observed that (1) when a large set of places visited was projected on the Euclidean plane, the empirical Hausdorff dimension of the set was always less than 2, and sometimes less than 1; (2) whereas traditional observation models assume reports at regular or exponentially spaced intervals, in real situations reports may be sporadic and bursty; (3) stochastic models based on memory-less random walks in space, as presented in \cite{ell} for example, are inadequate to represent a sequence of adaptive choices by an intelligent agent, such as a sequence of places or URLs visited. These observations are not new; see \cite{bar}, \cite{boy}, \cite{gon}, and references cited therein for related work. Assertion (1) is supported by inspection (www.mrlc.gov/nlcd2011.php) 
 of areas marked as {\em developed high intensity} in the U.S. National Land Cover Database \cite{hom}, and by analysis (courtesy of Dylan Molho) of
trip record data 
(obtained by the website fivethirtyeight.com
via Freedom of Information Law request on July
20, 2015) of Uber taxi pick-ups in New York City.

\subsection{A stochastic mobility simulator}
The first author has built and published a stochastic mobility simulator
called \texttt{FRACTALRABBIT} \cite{github}, with three tiers:

\begin{enumerate}
\item
    An \textit{Agoraphobic Point Process} generates a set $V$ of 
points in $\mathbb{R}^d$, whose limit is a random fractal, representing sites that could be visited.
\item
    A \textit{Retro-preferential Process} generates a trajectory $X$ through $V$, with strategic homing and self-reinforcing site fidelity as observed in human/animal behavior.
\item
    A \textit{Sporadic Reporting Process} to model time points $T$ at which the trajectory $X$ is observed, with bursts of reports and heavy tailed inter-event times.
\end{enumerate}
The present paper gives mathematical details about the first tier only. The second and third tiers will be described fully in later works.

\subsection{High level view of the agoraphobic point process}
\begin{figure}
\caption{\ \textit{Discrete Time Spatial Branching with Gaussian Displacements in $\R^2$ }:
Displacements follow a Gaussian law (\ref{eq:gauss}). Trees were truncated at 10,000 vertices.
Two  different exponents in (\ref{eq:rho}) were used,
giving random fractals with Hausdorff dimensions $10/9$ and $3/2$, respectively.
} \label{f:cgtree}
\begin{center}
\scalebox{0.25}{\includegraphics{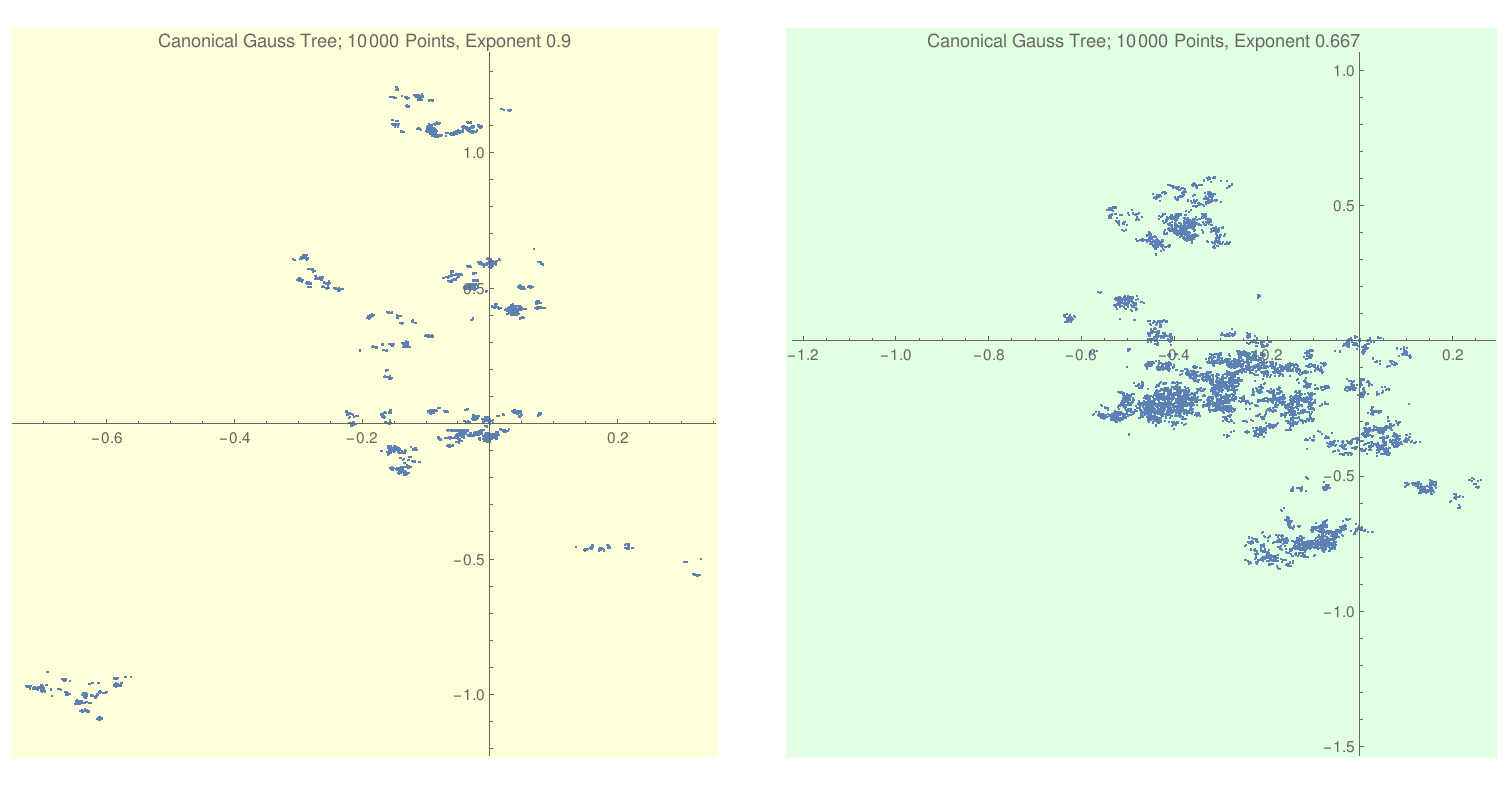}} 
\end{center}
\end{figure}

Our aim is to analyze the metric and dimensional properties of a 
certain random sets of points in $d$-dimensional euclidean space.  
These random sets are models for geographical clustering in which
new points are likely to be near to old points.  The models are not
time homogeneous.  Rather, as time progresses, the typical distance
between a new point and the set of existing points diminishes.
The parameters of the model will ensure that the set of points 
is finite at all times and countably infinite at time infinity;
the diameter remains bounded, whence there is a nontrivial limit
set; we are interested in how dense this set is, as measured by 
its Hausdorff dimension.  We begin with an informal description,
then proceed to a formal construction.

The process begins with a single point at the origin.  A density
$f$ on step sizes is fixed.  Each existing point $x$ gives rise 
to new points according to a spacetime Poisson process, where the 
total intensity at time $t$ (integrated over space) is proportional 
to $1/t$ and the displacement of the new point $y$ from the parent $x$ 
is isotropic with $|y-x|$ drawn from the density $f$ rescaled by 
$t^{-1/d}$.  We will see (Corollary \ref{cor:growth} below) that
the total number of points $\nt (t)$ grows like a random multiple
of $t^\rho$ where $\rho$ can be any positive constant and may be 
computed from $f$.  When the number of points reaches $n$, one has
$t \asymp n^{1/\rho}$ while the displacements of new points from old
ones is of order $t^{-1/d} \asymp n^{-1/(d \cdot \rho)}$.  Therefore, 
taking $\rho := \alpha / d$, our model is capable of emulating 
a purely spatial model in which the $(n+1)^{st}$ point is added 
by choosing uniformly among the existing points and adding an 
isotropic displacement with size distribution $f$ scaled by any 
exponent $n^{-1/\alpha}$, as described in Section \ref{s:compute}.
In particular, the model can be fitted to spatial data consisting
of a static collection of locations.

\section{Definitions and notation}

Because each new point is associated with a particular old point, 
it will be useful to represent the process as a tree, considering
the new point to be the offspring of the old point.  We begin with
notation for trees, then give a formal construction of the process
on an appropriate probability space.

\subsection{Notation for trees}

A rooted tree $T$ with vertices in a space $\vertex$ is a triple 
$(\rt, V, E)$ with $V$ a subset of $\vertex$ (the vertex set
of $T$), $\rt$ an element of $V$ (the root of $T$) and 
$E$ a subset of $V \times V$ (the set of oriented edges of $T$).
The set $E$ must satisfy the in-degree condition: there is
no edge $(v,\rt)$ in $E$ and for all $v \neq \rt$ in $V$ 
there is a unique $u \in V$ with $(u,v) \in E$.  This 
vertex $u$ is called the parent of $v$ and denoted $\parent (v)$.
The resulting graph must be connected, or in other words, $E$
must be well founded: for all $v \in V$ there exists a non-negative
integer $n$ such that iterating the parent map $n$ times returns
the root.  We denote this number by $|v|$; thus $\parent^{|v|} (v) = \rt$.
This formalism in sufficiently general to allow the number of children 
of a vertex to be finite, countably infinite or uncountable, however 
all of our trees will have vertex degrees that are at most countable.  

For $v, w$ vertices of $T$, write $v \to w$ to denote the relation 
that $v$ is the parent of $w$, and write $v \leq w$ to denote the
relation that $v$ is an ancestor of $w$ (the transitive closure of 
$v \to w$).  Let $v \wedge w$ denote the most recent common 
ancestor (the notation is consistent with the lattice meet 
in the ancestry partial order).  If $T$ is an infinite tree,
let $\partial T$ denote the set of infinite directed paths
$\gamma = (\rt, \gamma_1, \gamma_2, \ldots)$ from the root of $T$.  
For $\gamma \neq \gamma'$, let $\gamma \wedge \gamma'$ denote 
$\gamma_n$ where $n$ is the maximum of values such that 
$\gamma_n = \gamma_n'$. 
Recall some facts from \cite{LP-book} about the topology of
$\partial T$. The space $\partial T$ is topologized 
by a basis of clopen sets of the form 
$\partial T(v) := \{ \gamma : \gamma_n = v \}$.
The Borel $\sigma$-field with respect to this topology is denoted $\B$.
If $T$ has finite degrees then this topology  makes $\partial T$ 
a compact Hausdorff space.  
In general, one can compactify
by embedding in the space $\partial^+ T := \partial T \cup V(T)$ 
where a neighborhood basis of $v$ are the sets $T(v;F) := 
\{ \gamma : \gamma_n = v, \gamma_{n+1} \notin F \}$ for 
finite subsets $F$ of the children of $v$.  This makes
$\partial^+ T$ compact Hausdorff; it is still compact if
one removes the vertices of finite degree, as these are
isolated points in $\partial^+ T$.

The tree representing processes of interest here might naturally
be taken to have vertices in $\R^d$ or in $\R^d \times \R^+$,
representing locations or locations together with birth times.
However, because the graph structure is not random (each vertex 
has countably many children), it will be more convenient to 
fix a particular tree and define the random process as a random
function on the vertices of the canonical tree.  Accordingly,
we use the Ulam-Harris tree $\ulam$, whose 
vertex set is the set $\VU : = \bigcup_{n=1}^\infty \N^n$, with the
empty sequence $\rt = \emptyset$ as the root, with edges
from any sequence $\bfa = (a_1, \ldots , a_n)$ to any extension 
$\bfa \sqcup j = (a_1, \ldots , a_n , j)$ (see~\cite{louigi-ford}). 
Our construction will yield random maps $\tau : \VU \to \R^+$
and $\chi : \VU \to \R^d$ interpreted as birth times and locations.
Ultimately we will be interested in the range of $\chi$, or in
its closure or limit points.

\subsection{Probabilistic constructions}

Fix a positive integer $d$ which will be the spatial dimension.  
Fix positive real intensity parameters $\beta$ and $\theta$ and a 
spatial decay profile $f : \R^+ \to \R^+$.  One can without loss
of generality take $\beta = 1$, as the process, up to a linear 
time change, depends only on $\beta / \theta$. 
It is necessary that $f$ be integrable but not that it have 
total mass~1.  We assume $f$ is bounded and 
has finite moments of all orders. This holds, for example, if
$f$ has exponential tails, $f(x) < A e^{-Bx}$ for some $A,B > 0$.

Let $\m$ denote $d$-dimensional Lebesgue measure. 
Let $(\Omega , \F , \P)$ be a probability space on which is constructed 
a collection of IID Poisson point processes $\{ N_v : v \in \VU \}$ 
indexed by the canonical vertex space of the Ulam-Harris tree.  
The space of points for each point process is $\R^d \times \R^+$ 
and the common intensity is  
\begin{equation} \label{eq:gen pois}
F(y,t) := f \left [ (t/\beta)^{1/d} |y| \right ] 
   \frac{d\m (y) \, dt}{\theta} \, .
\end{equation}
Let
\begin{equation} \label{eq:c_d}
c_d := c_d (f) := \int_{\R^d} f(|x|) \, d\m (x)
\end{equation}
denote the total mass of the measure $F(y,1) \, d\m(y)$ in dimension $d$ in the case where 
$\beta / \theta = 1$.
By scaling, the total mass of $F(y,t) \, d\m(y)$ is $c_d / t$
when $\beta / \theta = 1$, and 
$c_d \beta/ (t \theta)$ in general.
Define
\begin{equation} \label{eq:rho}
\rho := c_d \frac{\beta}{\theta} \, .
\end{equation}

The random maps $\tau$ and $\chi$ are constructed recursively.
Begin with $\tau (\rt) := 1$ and $\chi (\rt)$ equal to the origin.
For the recursion, suppose that $\tau (v) = t$ and $\chi (v) = x$.  
The values of $\tau$ and $\chi$ on the children of $v$ are constructed 
from the Poisson process $N_v$ on $\R_d \times (t,\infty)$.  For
any $t' > t$, the total mass of $F$ on $\R^d \times (t,t')$
is finite, therefore $N_v (\R^d \times (t,t')) < \infty$ almost
surely for each $t' > t$, hence the points of the Poisson
process $N_v$ restricted to $\R^d \times (t,\infty)$ can be
enumerated in increasing order of the time coordinate: 
$(x_1, t_1), (x_2, t_2), \ldots$.  Now define 
\begin{eqnarray*}
\tau (v \sqcup j) & = & t_j \\
\chi (v \sqcup j) & = & x + x_j \, ,
\end{eqnarray*}
thereby extending the definitions of $\tau$ and $\chi$ to the
children of $v$ and completing the recursion.

The random map $\tau$ is equivalent to a well known branching 
random walk (BRW), called the Poisson Weighted Infinite Tree (PWIT),
introduced in~\cite{aldous-assignment}, with the name PWIT bestowed
in~\cite{aldous-zeta} (see also~\cite{aldous-steele-objective}).  
The PWIT is a BRW in the sense that there are variables 
$\{S(v) : v \in V\}$ with the collection 
$\{S(v) - S(\mbox{parent}(v))\}$ independent as parent$(v)$
varies over $V$.  
In the PWIT, the values $S(v \sqcup 1), S(v \sqcup 2) , \ldots$ 
at the children of $v$ are $S(v) + Y_1$, $S(v) + Y_2$, $\ldots$ 
where $\{ Y_n \}$ are the successive points of a rate~1 Poisson process.
We state the equivalence as a proposition.

\begin{proposition}[PWIT] \label{pr:PWIT}
For $v \in \VU$ and $n \geq 1$, let 
$$Y_n (v) := \log \frac{\tau (v \sqcup n)}{\tau (v)} \, .$$

\begin{enumerate}
\item[(i)] The vector $\rho \Y (v) = (0, \rho Y_1 (v) ,  \rho Y_2 (v) , \ldots)$ 
has IID mean~1 exponential increments, where $\rho = c_d \beta / \theta$.
\item[(ii)] The vectors $\Y (v)$ are independent as $v$ varies over $\VU$.
Consequently, $\Y$ are the increments of a PWIT. 
\item[(iii)] Let $\timesf$ be the $\sigma$-field generated by the PWIT,
that is, by all the vectors $\Y(v)$.  Conditional on $\timesf$, 
the increments $\Delta \chi (v) := \chi (v) - \chi (\parent (v))$ 
are independent with laws $F(x , \tau(v)) \, d\m(x)$ normalized 
to probability distributions.
\end{enumerate}
\end{proposition}

\noindent{\sc Proof:} Each $\Y (v)$ is constructed from the
Poisson process $N_v$, making~$(ii)$ automatic.  Projecting
$\R^d \times \R^+$ to $\R^+$, the intensity of 
$N_v (\R^d \times (t,\infty)$ is
$\frac{\rho}{t} \one_{t \geq \tau(v)} \, dt$.  
Making a time change from $t$
to $\log{t}$ results, after elementary change of variables, in a Poisson process of intensity
$\rho \one_{x \geq \log \tau (v)} \, dx$ and~$(i)$ follows
immediately.  For~$(iii)$, condition on a point of $N_v$ 
projecting to $[u , u+h]$ and let $h \downarrow 0$.
$\Cox$

Some examples of the spatial distributions we consider are as follows.
In all the examples, without loss of generality, $\beta$ is taken to be~1.
The volume of the unit ball is denoted $V_d := \pi^{d/2} / \Gamma (d/2 + 1)$.
\begin{example}[Exponential]
Let $f(x) = e^{-x}$.  Then 
\begin{equation} \label{eq:exp}
F(y,t) := \exp \left [ - t^{1/d} |y| \right ] 
   \frac{d\m (y) \, dt}{\theta} \, 
\end{equation}
and displacements of children from the parent will be exponential
with mean $t^{-1/d}$.  The constant $c_d$ is given by 
$$c_d = \int_0^\infty e^{-x} V_d \cdot d x^{d-1} \, dx = d! \, V_d \, .$$
\end{example}

\begin{example}[Gaussian]
Let $f(x) = e^{-x^2/2}$.  Then 
\begin{equation} \label{eq:gauss}
F(y,t) := \exp \left [ - \frac{|y|^2}{2 t^{-2/d}} \right ] 
   \frac{d\m (y) \, dt}{\theta} \, 
\end{equation}
and displacements of children from the parent will be centered
Gaussians with variance $t^{-2/d}$.  One obtains $c_d = \int 
e^{-|x|^2/2} \, d\m(x) = (2 \pi)^{d/2}$.
\end{example}

\begin{example}[Hard cutoff]
Let $f(x) = \one_{x \leq 1}$.  Then 
\begin{equation} \label{eq:hard}
F(y,t) := \one_{|y| \leq t^{-1/d}} \frac{d\m (y) \, dt}{\theta} \, 
\end{equation}
Here $c_d = V_d$, and displacements of children from the parent 
will be uniform within the ball of radius $t^{-1/d}$.
\end{example}

\subsection{Markov process representation}

To connect this construction to what was mentioned in the
introduction, we give an alternate description as a continuous time Markov process,
started at time $t = 1$.
For $t \geq 1$, let $\tree_t$ be the subtree of $\ulam$
induced on those vertices $v$ satisfying $\tau (v) \leq t$.
Without ambiguity, the restrictions of $\tau$ and $\chi$ to $\tree_t$
are also denoted $\tau$ and $\chi$.  Let $\F_t$ denote the 
$\sigma$-field generated by the restrictions of all 
Poisson processes $N_v$ to $\R^d \times [1,t]$.  Let $\pace_0$ 
be the set of finite trees with vertices in $\R^d \times \R^+$. 
Without loss of generality, we identify $\tree_t$ with the
tree in $\pace_0$ whose vertex set is the range of $(\chi, \tau)$
restricted to $\tree_t$ with edges such that $(\chi,\tau)$ 
is a graph isomorphism.  Let $\zerotree$ denote the rooted tree 
whose unique vertex is the pair $({\bf 0} , 1)$.
The following result is self-evident.
\begin{proposition} \label{pr:Markov}

\begin{enumerate}
\item[(i)]
The process $\{ \tree_t : t \geq 1 \}$, together with the maps $\tau$
and $\chi$, is Markov with respect to the filtration $\{ \F_t \}$.
It is a time-inhomogeneous pure jump process on $\pace_0$ with 
initial state $\zerotree$ whose jumps at time $t$ add vertices $(x,t)$ 
with kernel $F(x-v , t) \, d\m(x)$.
\item[(ii)]
Projecting $\R^d \times \R^+$ to the last coordinate
yields a Markov process which is a BRW in the sense that any
vertices alive at time $t$ independently give birth to new 
vertices, marked with the timestamp $t$, at rate $\rho/t$. 

\item[(iii)]
The size $\nt_t$ of $\tree_t$ is a pure birth chain with
birthrate $\rho \nt_t / t$.
\end{enumerate}
\end{proposition}

The following consequence makes precise the heuristic discussion 
in the introduction
\begin{corollary}[growth exponent] \label{cor:growth}
With probability~1, the limit
$$W := \lim_{t \to \infty} t^{-\rho} \nt (t)$$
exists and is nonzero.
\end{corollary}

\noindent{\sc Proof:}
Let $L_t := \log \nt (t)$.  The jump rate $\rho \nt (t) / t$ for 
$\{ \nt  (t) \}$ implies a jump rate $\rho \exp(L(t)) / t$ for
$\{ L(t) \}$; the increment is $\log [1 + \exp (-L(t))]$.
Then $L_t$ is an increasing pure jump process with infinitesimal
drift and variance give by 
\begin{eqnarray*}
\mu_t & = & \frac{\rho}{t} e^{L_t} \log (1 + e^{-L_t}) \\
\sigma^2_t & \leq & \frac{\rho}{t} e^{L_t} e^{-2L_t} \, .
\end{eqnarray*}
We first establish a crude bound: for some $\alpha > 0$,
\begin{equation} \label{eq:crude}
W' := \sup_t t^{\alpha} e^{-L_t} < \infty \, .
\end{equation}

Indeed, $\mu_t > \rho / (2t)$ for any $t \geq 1$, whence for any 
$\alpha < \rho/2$, for sufficiently large $t$ we have $L_t > \alpha \log t$ 
and~\eqref{eq:crude} follows.  From~\eqref{eq:crude}, 
it follows that if $\tau_c := \inf \{ t : t^\alpha e^{-L_t} \geq c \}$
then the events $\{ \tau_c < \infty \}$ converge down to zero as 
$c \to \infty$.  Taking Taylor series, it follows on $\{ \tau_c = \infty \}$
that 
\begin{eqnarray}
\mu_t & = & \frac{\rho}{t} + O(t^{-1-\alpha}) \label{eq:mu} \\
\sigma_t^2 & \leq & c \rho t^{-1-\alpha} \, . \label{eq:sigma}
\end{eqnarray}
Let $M_t := L_t - \int_1^t \mu_s ds$.  Then $M_{t \wedge \tau_c}$ is a
martingale in $L^2$ for every $c$.  It follows from~\eqref{eq:sigma} 
and the fact that $\tau_c = \infty$ for $c > W'$
that the almost sure
limit $M_\infty := \lim_{t \to \infty} M_t$ exists for all $t$.  
Then from~\eqref{eq:mu} we get $\int_1^\infty (\mu_t - \rho / t) \, dt$
converges almost surely to a finite random variable $\Delta$, whence
$$L_t - \rho \log t \to M_\infty + \Delta$$ 
proving the proposition with $W = \exp(M_\infty + \Delta)$.
$\Cox$

\subsection{Metric definitions}

Let $\delta$ be a metric on $\partial^+ \ulam$ compatible with the
topology generated by the clopen sets $\partial T(v)$.  One such
metric is the metric $\delta$ defined by 
$\delta (\gamma , \gamma') := \tau(\gamma \wedge \gamma')^{-1/d}$.
Finiteness of the set $\{ v : \tau (v) < t \}$ for every $t$,
which occurs with probability one, ensures that this random metric
is indeed a metric and generates the correct topology, which is compact.
Any metric such as $\delta$ which depends only on $\gamma \wedge \gamma'$
is an {\em ultrametric}, meaning that among the three values
$\delta(\gamma , \gamma')$, $\delta(\gamma , \gamma'')$ and 
$\delta(\gamma', \gamma'')$, the greatest two are always equal.

We recall the definition of the Hausdorff dimension of a closed
subset $S$ of a compact metric space. 
\begin{defn}[Hausdorff dimension] \label{def:dim}
The $\alpha$-dimensional {\em Hausdorff content} $H_\alpha (S)$ is
the infimum value of $\sum_{j=1}^\infty \diam (B)^\alpha$ over 
collections of sets $\{ B(x_j , r_j) : j \geq 1 \}$ covering $S$.
The {\em Hausdorff dimension} $\dim (S)$ is the supremum of $\alpha$ 
for which $H_\alpha (S) = \infty$ and also the infimum of $\alpha$ 
for which $H_\alpha (S) = 0$.  
\end{defn}

\begin{remarks}[Hausdorff measure]
Finer information can be obtained by using gauge functions other
than the $\alpha$ power.  We will not need these here.  When 
$\alpha = \dim (S)$, it can happen that $H_\alpha (S)$ is zero, or
finite nonzero.  The $\alpha$-dimensional 
{\em Hausdorff measure} of $S$, defined by the increasing limit
as $\ee \to 0$ of $H_\alpha (S)$ when the infimum over covers 
is restricted to balls of diameter less than $\ee$, may be zero, finite
or infinite, according to how the content behaves; while Hausdorff measure is
more complicated than Hausdorff content by one limit, 
it has the
benefit when finite of being a measure, that is, additive over 
disjoint sets.
\end{remarks}

We recall from~\cite{LP-book} an equivalent definition of 
Hausdorff dimension.  First they show (e.g. as a consequence 
of~\cite[Theorem~15.7]{LP-book}) that in the infimum 
over covers, it suffices to consider only covers by sets
of the form $\partial T(v)$. Indeed for ultrametrics, as here,
one can assume without loss of generality that every
element of a finite cover is of the form $\partial T(v)$.

Next, consider finite directed 
flows on $T$ defined by placing positive real numbers on
the directed edges so that inflow equals outflow at every 
vertex, except at the root where total outflow is finite.
Such flows are obviously in bijective correspondence with
finite Borel measures on $\partial T$.  Given $\alpha$,
consider the class of {\em admissible} flows, constrained
never to exceed $\delta(v)^\alpha$ at the vertex $v$; 
here $\delta(v)$ is short for $\delta (\gamma , \gamma')$ 
for any/every pair with $\gamma \wedge \gamma' = v$.

By the max-flow min-cut theorem for countable directed networks
~\cite[Theorem~3.1]{LP-book}), the supremum of volume for
admissible flows has magnitude equal to the infimum of the
sum of constraints over cutsets, and there is a flow 
attaining this supremum.  But the sum of constraints over 
cutsets is precisely $H_\alpha (\partial T)$.  
We have therefore shown:

\begin{lemma}[Hausdorff dimension in terms of flows] \label{lem:flows}
Suppose $\delta$ is the metric for which $\delta (\gamma , \gamma')$ is $\tau (\gamma \wedge \gamma')^{-1/d}$.
The Hausdorff dimension of $\partial T$ is the critical $\alpha$
for whether it is possible to have a nonzero flow constrained
to be at most $\delta(v)^\alpha$ at each vertex $v$.
$\Cox$
\end{lemma}

One further equivalent definition will be useful, in terms 
of capacity.  Define the $\alpha$-dimensional energy of 
a measure $\mu$ on $\partial T$ by
\begin{equation} \label{eq:energy}
\cE_\alpha (\mu) := \int \int \delta(\gamma,\gamma')^{-\alpha} 
   \, d\mu (\gamma) \, d\mu (\gamma')  \, .
\end{equation}
The following criterion for bounding dimension from below
based on existence of measures of finite energy may
be found in~\cite{LP-book} or~\cite{bishop-peres}.
\begin{lemma}[Capacity definition of Hausdorff dimension] \label{lem:cap}
Let $A$ be a closed subset of $\partial T$.  If there is a 
nonzero measure $\mu$ on a set $A$ with $\cE_\alpha (\mu) < \infty$ 
then $\dim (A) \geq \alpha$.
\end{lemma}

\noindent{\sc Proof:} 
By definition, we take
$\delta(\gamma,\gamma) = \infty$, hence any finite energy measure must be non-atomic.  Thus 
\(\mu \times \mu\) 
is supported on pairs \((\gamma , \gamma')\) with \(\gamma \neq \gamma'\) .  Therefore
some small multiple of $\mu$, interpreted
as a flow, is admissible for the constraints $\{ \delta(v)^\alpha \}$, 
and thus witnesses the criterion in Lemma~\ref{lem:flows}.  
In fact, therefore, not only is $\alpha$ 
no more than the critical dimension, 
but if $\alpha$ is the critical dimension, then 
the Hausdorff content (and hence Haudorf measure) in the critical dimension is positive. 
$\Cox$

\section{Dimension of $\tree_\infty$ in metric $\delta$}

\begin{thm} \label{th:dim 1}
With respect to the metric $\delta$, the dimension of 
$\partial \ulam$ is $d \cdot \rho$.
\end{thm}

\noindent{\sc Proof:}  
The set of constraints $\{ \delta (v)^\alpha : v \in \VU \}$
form a {\em Galton-Watson} network in the sense 
of~\cite[Section~5.9]{LP-book}; a slight generalization 
is required to allow infinitely many children. 
The weights $\{ A_j \}$ in~\cite{LP-book} are the constraints
$(t_1^{-\alpha/d} , t_2^{-\alpha / d} , \ldots)$.  Theorem~5.35
of~\cite{LP-book} states an intuitively obvious result, namely 
water flows if $\E \sum A_i > 1$ and not if $\E \sum A_i < 1$
(the case $\E \sum A_i = 1$ is not settled by that result).
In particular, the critical dimension for water flow is the 
$\alpha$ for which $\E \sum A_i = 1$.  

Recall that the birth times for children of the root is a Poisson
point process with intensity $\rho \, dt / t$.  We may therefore
compute
\begin{equation} \label{eq:sum}
\E \sum_i A_i = \int_1^\infty \rho t^{-\alpha/d} \frac{dt}{t}
   = \rho \int_1^\infty t^{-\alpha/d - 1} = \rho \frac{d}{\alpha} \, .
\end{equation}
We see that $\alpha = \rho d$ is critical for $\E \sum_i A_i \geq 1$,
hence for water flow with capacities $t^{-\alpha / d}$.
By Lemma~\ref{lem:flows}, the dimension of $\partial \ulam$ 
with metric $\delta$ is thus shown to be $\rho \cdot d$.
$\Cox$

Before continuing, we give a second proof of the lower dimension
bound in Theorem~\ref{th:dim 1}.  This will provide a framework we
will use to prove Theorem~\ref{th:HD} below.  
First, we construct a martingale for the Galton-Watson network
from the proof of Theorem~\ref{th:dim 1} analogous to the normalized
size martingale for Galton-Watson trees.
Define $Z(v) := \tau (v)^{-\rho}$.  We have seen in~\eqref{eq:sum} that 
$\E \sum_{|v|=1} Z(v) = 1$, where
 $|v|$ denotes the depth of $v$ in the tree.
A similar computation shows that 
\begin{equation} \label{eq:moment}
\E \sum_{|v|=1} \tau (v)^{-\lambda} 
   = \rho \int_1^\infty t^{-\lambda - 1} \, dt = \rho/\lambda \, ,
\end{equation}
and, by induction, that
$$\E \sum_{|v|=n} \tau (v)^{-\lambda} = 
   \left ( \frac{\rho}{\lambda} \right )^n \, .$$
Applying this with $\lambda = 2 \rho$ shows that 
$\Var \sum_{|v|=1} Z(v) = \E \sum_{|v|=1} \tau (v)^{2\rho} = 1/2$.  
Define $W_n := \sum_{|v|=n} Z(v)$.  Letting $z \geq y$
denote the relation holding when $y$ is an ancestor of $z$, generalize this by defining
$$W_n (y) := \sum_{z \geq y : |z| = |y| + n} \frac{Z (z)}{Z(y)} \, .$$
For $y$ varying over vertices at a fixed depth $k$, the variables
$W_n (y)$ are IID with the same distribution as $W_n$.  By
convention, if $n=0$, we take $W_n(y) = 1$.  From the definitions 
and induction, it is easy to see that for any $1 \leq k \leq n$,
\begin{equation} \label{eq:recur}
W_n (x) = \sum_{y \geq x, |y| = k+|x|} Z(y) W_{n-k} (y) \, ,
\end{equation}
summing over descendants $y$ of $x$.
It follows from this that $\E (W_n \| \F_{n+1}) = W_{n+1}$ where
$\F_n$ is the $\sigma$-field of information in the first $n$
levels of the tree.  Thus $\{ W_n \}$ is a martingale.
These variables are square integrable, as shown by the following 
computation.
\begin{eqnarray*}
\E W_n^2 & = & \E \left ( \sum_{|x|=1} Z(x) W_{n-1}(x) \right )^2 \\[1ex]
& = & \sum_{|x|=1} \E Z(x)^2 W_{n-1}(x)^2 
   + \sum_{|x|=|y|=1, x \neq y} \E Z(x) Z(y) W_{n-1}(x) W_{n-1}(y) \\[1ex]
& = & \E W_{n-1}^2 \E \sum_{|x|=1} Z(x)^2 +
    \sum_{|x|=|y|=1, x \neq y} \E Z(x) Z(y) \\[1ex]
& = & \E (\sum_{|x|=1} Z(x))^2 + \left ( 
   \E \sum_{|x|=1} Z(x)^2 \right ) \left ( \E W_{n-1}^2 - 1 \right ) \\[1ex]
& = & \frac{3}{2} + \frac{1}{2} \left ( \E W_{n-1}^2 - 1 \right ) \, .
\end{eqnarray*}
Inductively, $\E W_{n-1}^2 = 2 - 2^{-n}$, hence 
$\lim_{n \to \infty} \E W_n^2 = 2$.  
Therefore the martingale $\{ W_n \}$ is square integrable,
and so it converges almost surely and in $L^2$ to some random variable $W$ with mean~1 and $\E W^2 = 2$.
Similarly, $W(v) := \lim_{n \to \infty} W_n (v)$ 
defines a limiting random variable associated with each vertex.

We summarize the foregoing computation in a proposition.

\begin{proposition}[Limit uniform measure] \label{pr:mart}
For each $v \in \VU$ there is a martingale $\{ W_n (v) \}$ converging 
almost surely and in $L^2$ to a limit $W(v)$ with $\E W(v) = 1$
and $\E W(v)^2 = 2$.  Define
\begin{equation} \label{eq:LU}
\mu (V) := Z(v) W(v) \, .
\end{equation}
It follows from~\eqref{eq:recur} that $\mu$ is a measure on $\partial \ulam$.
We call this the $d \cdot \rho$-\textbf{dimensional limit uniform measure}.
$\Cox$
\end{proposition}

\noindent{\sc Second proof of lower dimension bound:}
By Lemma~\ref{lem:cap} it suffices to show for every $a < d \cdot 
\rho$ that $\cE_a (\mu) < \infty$.  Bound $\cE_a$ from above by
\begin{eqnarray*}
L & := & \sum_v \mu (\partial T(v))^2 \tau (v)^{a/d} \\[1ex]
& \geq & \sum_v (\mu \times \mu) \, 
   \{ (\gamma , \gamma') : \gamma \wedge \gamma' = v \} \; \tau (v)^{a/d} \\[1ex]
& = & \cE_a (\mu) \, .
\end{eqnarray*}
Taking expectations and breaking down the sum according to depth, 

\begin{eqnarray}
\E L & = & \sum_n \E \sum_{|v| = n} \mu (\partial T(v))^2 \tau (v)^{a/d} 
   \nonumber \\[1ex]
& = & \sum_n \E \sum_{|v|=n} \tau (v)^{-2 \rho} W(v)^2 \tau (v)^{a/d} 
   \nonumber \\[1ex]
& = & 2 \sum_n \E \sum_{|v|=n} \tau (v)^{-\rho - \ee} \label{eq:tree energy}
\end{eqnarray}
where $\ee := (d \cdot \rho - a)/d > 0$ by assumption.  In~\eqref{eq:moment}
we saw that the inner sum is $(\rho / (\rho + \ee))^n$.  Hence the outer sum 
is a geometric series summing to $1 + \rho/\ee < \infty$.  This completes
the proof that $\cE_a (\mu) < \infty$ for $a < d \cdot \rho$ and hence that
$\dim (\partial \tree_\infty) = d \cdot \rho$.
$\Cox$

\section{Dimension of the euclidean set}

We turn now to the question of the dimension of the 
euclidean set, that is, the closure of the range of $\chi$.  
It is least complicated still to work on $\partial \ulam$,
identifying $\gamma \in \partial \ulam$ with the point 
$\chi (\gamma) := \lim_{n \to \infty} \chi (\gamma_n)$.  
We do this by lifting the euclidean metric to 
the metric $\delta'$ on $\partial \ulam$ defined by
$$\delta' (\gamma , \gamma') := |\chi (\gamma) - \chi (\gamma')| \, .$$
This is not an ultrametric, therefore using Lemma~\ref{lem:flows} 
can give only an upper dimension bound.  This is not a problem because
we are going to use Lemma~\ref{lem:cap} for the lower bound.
Our result is that mapping $\partial \ulam$ into euclidean space 
does not reduce the dimension unless it must, due to the dimension
of $\partial \ulam$ being greater than $d$.
\begin{thm}[euclidean dimension] \label{th:HD}
The Hausdorff dimension of $\partial \ulam$ with
respect to $\delta'$ is $d \cdot \min \{ 1 , \rho \}$.
\end{thm}

\noindent{\sc Proof of Lower dimension bound:}
Fix $a < \min \{ d , d \cdot \rho \}$.  We need to show that the
Hausdorff dimension of $\partial \tree_\infty$ in the metric $\delta'$
is at least $a$.  Let $\mu$ denote the $d \cdot \rho$-dimensional
limit uniform measure.  We show that the expected $a$-dimensional 
energy with respect to $\delta'$, $\cE_a' (\mu)$ is finite.  Finite
expectation implies almost sure finiteness, which by Lemma~\ref{lem:cap}
implies $\dim(\partial T)$ in the metric $\delta'$ is almost surely
at least $a$.  We do this in two steps.  Given $\gamma , \gamma' \in
\partial T$, if $\gamma_j = \gamma_j'$ for $j \leq n$ but $\gamma_{n+1}
\neq \gamma_{n+1}'$, denote
$$\delta'' (\gamma , \gamma') := 
   \min \{ \tau (\gamma_{n+1}) , \tau (\gamma_{n+1}') \}^{-1/d} \, .$$
By comparison, $\delta (\gamma , \gamma') = \tau (\gamma_n)^{-1/d}$, 
therefore $\delta'' < \delta$, because the $-1/d$ power is of the
first time occurring in only one of the two branches, rather than 
of the last common time.  

\begin{lemma} \label{lem:1a}
The measure $\mu$ has finite $a$-dimensional energy in metric $\delta''$:
$$\E \cE_a'' (\mu) < \infty \, .$$
\end{lemma}

Assuming this for the moment, the lower dimension bound is proved as follows.  
Let $\G$ be the $\sigma$-field defined generated by the time variables 
$\{ \tau (v) : v \in \VU \}$.  Given $\gamma$ and $\gamma'$, let
$v = \gamma \wedge \gamma'$, let $n = |v|$, and let $i$ and $j$ be
the distinct positive integers for which $\gamma_{n+1} = v \sqcup i$
and $\gamma_{n+1}' = v \sqcup j$.  Write $\gamma \prec \gamma'$
if $i < j$.  For any measure, and in particular for $\mu$,
\begin{equation} \label{eq:''}
\cE_a'' (\mu) = 2 \sum_v \sum_{i < j} 
   \mu (\partial T(v \sqcup i)) \mu (\partial T(v \sqcup j)) \tau (v \sqcup i)^{a/d} \, .
\end{equation}
By comparison,
\begin{equation} \label{eq:'}
\cE_a' (\mu) = 2 \sum_v \sum_{i < j} 
   \int |\chi(\gamma) - \chi (\gamma')|^{-a} \, 
   d\mu |_{T(v \sqcup i)} (\gamma) \, d\mu |_{T(v \sqcup j)} (\gamma') \, .
\end{equation}
We will show that the expectation given $\G$ of~\eqref{eq:'} may 
be bounded term by term by the summands in~\eqref{eq:''}.  
In fact we claim there is a constant $K$ for which
\begin{equation} \label{eq:K}
\E \left ( \left. |\chi (\gamma) - \chi (\gamma')|^{-a} 
   \right | \| \G \right ) \leq K \tau (v \sqcup i)^{a/d} \, .
\end{equation}
To see this, let $\G' \supseteq \G$ be the $\sigma$-field
generated by all the variables $\tau (v)$ and all the variables
$\Delta \chi (v)$ except for $\Delta \chi (v \sqcup i)$.
It suffices to show~\eqref{eq:K} with $\G'$ in place of $\G$.
Let $x := \chi (\gamma) - \chi (\gamma') - \Delta \chi (v \sqcup i)$, and write
$$\E \left ( \left. |\chi (\gamma) - \chi (\gamma')|^{-a} 
   \right | \| \G' \right ) = \E |x + \tau(v \sqcup i)^{-1/d} \theta |^{-a}$$
where $\theta$ is a random variable with density $c_d^{-1} f(|x|) \, d\m(x)$; 
see Proposition~\ref{pr:PWIT}, part~$(iii)$.
Evidently $\tau (v \sqcup i)$ is measurable with respect to 
$\G \subseteq \G'$. 
The quantity $x$ is the sum of all displacements of children from parents along $\gamma $, 
minus the same sum from $\gamma'$, except that one summand from $\gamma$ gets deleted, namely the displacement of $v \sqcup i$ from $v$.  
Therefore, $x$ does not require this displacement to be evaluated and is measurable with respect to $\G'$ also.  
Abbreviate $\tau (v \sqcup i)$ to $\tau$,
and noting that $\theta$ is the only random quantity on the left side,
 we have a tail bound
\begin{eqnarray*}
\P (|x + \tau^{-1/d} \theta|^{-a} > \lambda) & = &
   \P (|x + \tau^{-1/d} \theta| < \lambda ^{-1/a}) \\[1ex]
   & = & \P (|\tau^{1/d} x + \theta| < \tau^{1/d} \lambda^{-1/a}) \\[1ex]
   & = & \P (|\tau^{1/d} x + \theta|^d < \tau \lambda^{-d/a}) \\[1ex]
   & \leq & 1 \wedge C \tau \lambda^{-d/a}
\end{eqnarray*}
because the density $f$ is bounded.
where the constant $C$ depends only on $f$.  By assumption
$-d/a < -1$.  Therefore, letting $\lambda_0 := (C \tau)^{a/d}$,
we obtain~\eqref{eq:K} via
\begin{eqnarray*}
\E |x + \tau^{-1/d} \theta|^{-a} & \leq &
   \int_0^\infty 1 \wedge C \tau \lambda^{-d/a} \, d\lambda \\[1ex]
& = & \lambda_0 + \int_{\lambda_0}^\infty C \tau \lambda^{-d/a} \, d\lambda 
   \\[1ex]
& = & C' \lambda_0 = K \tau^{a/d} \, .
\end{eqnarray*}

Integrating the left side of~\eqref{eq:K} against $\mu$ restricted to
$T(v \sqcup i)$ in the variable $\gamma$ and $\mu$ restricted to 
$T(v \sqcup j)$ in the variable $\gamma'$ then shows that
$$\E \left ( \left. \int |\chi (\gamma) - \chi (\gamma')|^{-a} 
   \, d\mu |_{T(v \sqcup i)} (\gamma) 
   \, d\mu |_{T(v \sqcup j)} (\gamma') \right | \G \right ] 
   \leq K \mu (v \sqcup i) \mu (v \sqcup j) \tau (v \sqcup i)^{-a} \, ,$$
bounding~\eqref{eq:'} above, term by term, by $K$ times~\eqref{eq:''}.
Thus, 
$$\cE_a' (\mu) = \E \left [ \E_{~} (\cE_a' (\mu) \, |^{~} \G) \right ] \leq
   K \cE_a'' (\mu) < \infty$$
by Lemma~\ref{lem:1a}, finishing the lower dimension bound.  
It remains to prove the lemma.

\noindent{\sc Proof of Lemma}~\ref{lem:1a}: For $\mu$ or any measure,
we may break up $\cE_a'' (\mu)$ as
$$\frac{1}{2} \cE_a'' (\mu) = \sum_n \sum_{|v| = n} \sum_{i < j}
   \mu (v \sqcup i) \mu (v \sqcup j) \tau (v \sqcup i)^{a/d} \, .$$
Let $\ee$ denote the positive quantity $\rho - a/d$ and let
$\theta_i = \theta_i (v)$ denote $\tau (v \sqcup i) / \tau (v)$.
Plugging in $\mu (v) = \tau(v)^{-\rho} W(v)$ then gives
\begin{eqnarray*}
\cE_a'' (\mu) & = & \sum_n \sum_{|v|=n} \tau(v)^{-2\rho} \tau(v)^{a/d}
   \sum_i \left [ 
      W_i \theta_i (v)^{-\rho + a/d} 
      \sum_{j > i} W_j \theta_i (v)^{-\rho}
   \right ] \\[2ex]
& = & \sum_n \sum_{|v|=n} \tau(v)^{-\rho - \ee} 
   \sum_i \left [ W_i \theta_i^{-\ee} 
      \sum_{j > i} W_j \theta_j (v)^{-\rho} \right ] \, .
\end{eqnarray*}
For each $v$, the quantities $\tau (v), \theta_i (v), \theta_j (v),
W(v \sqcup i)=:W_i$ and $W(v \sqcup j)$ are independent.  
Also, $\theta_i (v)$ is always less than~1 and $\E W_i = 1$ for all $i$.  
Therefore, reversing the inner summations,
\begin{eqnarray*}
\E \cE_a'' (\mu) & < & \sum_n \sum_{|v|=n} \E \tau(v)^{-\rho - \ee} 
   \E \sum_j \left [ W (v \sqcup j) 
   \theta_j (v)^{-\rho} \sum_{i < j} W_i \right ] \\[1ex]
& = & \sum_n \sum_{|v|=n} \E \tau(v)^{-\rho - \ee} 
   \sum_j (j-1) \E \theta_j (v)^{-\rho} \, .
\end{eqnarray*}
From $\E \theta_j^{-\rho} = (\E \theta_1^{-\rho})^j$ and $\theta_1 < 1$, 
it follows that $\sum_j (j-1) \E \theta_j^{-\rho}$ is equal to the 
constant $M := (\E \theta_1^{-\rho} / (1 - \E \theta_1^{-\rho}))^2$.
Therefore,
$$\E \cE_a'' (\mu) < M \E \sum_n \sum_{|v|=n} \tau (v)^{-\rho-\ee}$$
which was seen to be finite in~\eqref{eq:tree energy}, finishing
the lemma, hence the lower dimension bound.
$\Cox$

\noindent{\sc Proof of upper dimension bound.}  We may assume that
$\rho < 1$ because otherwise the upper dimension bound of $d$ is trivial.
Fix $a > \rho \cdot d$.
We need to show that the $a$-Hausdorff content of $\tree_\infty$ is zero.  
We use the covers $\{ \partial \tree (v) : |v| = n \}$.  
We will show 
that the expectations go to zero
$$\E \sum_{|x| = n} \diam (\partial \tree (v))^a < \infty \to 0$$
which implies with probability~1 the lim inf of the sums goes to zero, 
hence $H_a (\partial \tree_\infty) = 0$.  Because 
$\diam (\partial \tree (v)) / \tau (v)$ is
independent of $\tau (v)$ and distributed as $\diam (\partial \tree_\infty)$ 
we may write
\begin{eqnarray*}
\E \sum_{|x| = n} \diam (\partial \tree (v))^a & = & \E \sum_{|x| = n} 
   \tau (v)^a \left ( \frac{\diam (\partial \tree (v))}{\tau (v)} \right )^a \\[1ex]
& = & \left ( \E \sum_{|x| = n} \tau (v)^a \right )
   \E \diam (\partial \tree_\infty)^a \\[1ex]
& = & \left ( \frac{d \cdot \rho}{a} \right )^n 
   \E \diam (\partial \tree_\infty)^a \, .
\end{eqnarray*}

Because $d \cdot \rho / a < 1$, 
it suffices to show that $\diam (\ulam)$ has finite $a$-moment.
We show in fact all moments of $\diam (\ulam)$ are finite.  
We use a cheap upper bound, namely 
\begin{equation} \label{eq:chi}
\diam (\ulam) \leq 2 \sup_{\gamma \in \partial \ulam} 
   \sum_{n=1}^\infty |\chi (\gamma_n) - \chi (\gamma_{n-1})|
   \leq 2 \sum_{n=1}^\infty \sup_{|v| = n} |\chi (v) - \chi (\parent (v))| 
   \, .
\end{equation}
Let $M_n$ denote the supremum inside the sum.  We need tail bounds 
on $M_n$.  Recall from Proposition~\ref{pr:PWIT} that the values 
$S(v) := \log \tau (v)$ form a PWIT scaled by $1/\rho$.  By construction, 
conditional on $\timesf$ (that is all the values $\tau (v)$)
the displacements $|\chi (v) - \chi (\parent(v))|$ are independent, 
the one at each vertex $v$ chosen from $f$ scaled by $\exp (- S(v)/d)$.  
Thus the collection $\{ \chi (v) - \chi (\parent (v)) \}$ is just
$\{ Y(v) := \theta (v) \exp (-S(v)/d) \}$ where $\{ \theta (v) \}$ are 
IID picks from the density described by $f$.  

Lemma~3.2 of~\cite{louigi-ford} says that
\begin{equation} \label{eq:ABFo2013}
G_n (x) := \E \# \left \{ v : |v| = n \mbox{ and } PWIT(v) \leq x \right \} 
   = \frac{x^n}{n!} \, .
\end{equation}
The expected sum of the $m$-powers at level $n$ may be written as
$$E_n^m := \E \sum_{|v|=n} \{ Y(v)^m : |v| = n \}
   = \left [ \E \theta(v)^m \right ] \;  
   \E \sum_{|v| = n} e^{- (m/d) S(v)}$$
where $S(v)$ are $1/\rho$ times the values of the PWIT at $v$.
Expressing the last expectation as an integral against the 
density $G_n(x)$ of values of the PWIT at level $n$ and 
using~\eqref{eq:ABFo2013} then gives 
\begin{eqnarray*}
E_n^m & = & \left [ \E \theta(v)^m \right ] 
   \int_0^\infty e^{-\frac{a}{d\rho} x} \, dG_n (x) \\[1ex]
& = & \left [ \E \theta(v)^m \right ] 
   \int_0^\infty e^{-\frac{a}{d \cdot \rho} x} \, 
   \frac{x^{n-1} \, dx}{(n-1)!} \\[1ex]
& = & \left [ \E \theta(v)^m \right ] 
   \left ( \frac{a}{d \rho} \right )^n \, .
\end{eqnarray*}
Let $c_{m,f}$ be the $m$th moment of a random variable with density
$f$, normalized, which we have assumed to be finite for all $m$.
Markov's inequality then yields, for any $m$,
\begin{equation} \label{eq:markov}
\P (M_n \geq \lambda) \leq \lambda^{-m} E_n^m
   = \lambda^{-m} c_{m,f} \left ( \frac{m}{d \cdot \rho} \right )^{-n} \, .
\end{equation}

Fixing any $m > d \cdot \rho$, we may choose $\ee > 0$ such that 
$(1 + \ee)^m d \rho / m < 1$.  Let $y_n := (1 + \ee)^{-n}$ and 
apply~\eqref{eq:markov} to see that
\begin{equation} \label{eq:tail}
\P (M_n \geq \lambda \ee (1 + \ee)^{-n}) \leq
   c_{m,f} \ee^{-m} \lambda^{-m} 
   \left ( \frac{(1+\ee)^m}{m / (d \cdot \rho)} \right )^n \, .
\end{equation}
Recall that $\sum_n M_n$ is an upper bound on $\diam (\ulam)$,
which we are trying to show has all moments finite.
This follows if $\P (\sum_n M_n > \lambda) = O(\lambda^{-c})$
for every $c > 0$ as $\lambda \to \infty$.  Because $\sum_{n=1}^\infty
\ee (1+\ee)^{-n} = 1$, if $\sum_n M_n > \lambda$ then $M_n \geq
\lambda \ee (1+\ee)^{-n}$ for some $n \geq 1$.  Thus for any $m$, 
using~\eqref{eq:tail},
\begin{eqnarray*}
\P \left ( \sum_{n=1}^\infty M_n \geq \lambda \right ) & \leq &
   \sum_{n=1}^\infty \P \left ( M_n \geq \ee (1 + \ee)^{-n} \lambda
   \right ) \\[1ex]
& \leq & C \lambda^{-m} 
\end{eqnarray*}
where $\disp C := c_{m,f} \ee^{-m} \frac{(1+\ee)^m}{m/(d \cdot \rho) 
- (1+\ee)^m}$.  This finishes the proof that $\diam (\ulam)$ has 
finite moments, hence the proof that the $a$-dimensional Hausdorff 
content of $\partial \tree_\infty$ is zero for any $a > d\rho$,
establishing the upper dimension bound.
$\Cox$

As a corollary, we get the convergence of $\chi (\gamma_n)$ along
every path $\gamma$.

\begin{corollary} \label{pr:converges}
With probability~1, for every $\gamma \in \partial \tree_\infty$, 
the sequence $\{ \gamma_n \}$ converges.
\end{corollary}

\noindent{\sc Proof:} We have seen in~\eqref{eq:markov} 
that $\P (M_n \geq \ee (1+\ee)^{-n} \lambda)$ is summable for
some $\ee > 0$.  For any $\gamma$, the triangle inequality gives
$|\chi (\gamma_n) - \chi(\gamma_m)| \leq \sum_{j=m+1}^n M_j$.  Together 
with Borel-Cantelli, this implies that $\{ \chi (\gamma_n) : n \geq 1 \}$
is a Cauchy sequence simultaneously for all $\gamma$.
$\Cox$

\section{Computational Implementation} \label{s:compute}
\subsection{The Agoraphobic Point Process in the Unit Ball} \label{s:svmp}
The preceding theory described a continuous time branching
process embedded in $\R^d$. This section describes two finite
discrete time constructions used in \cite{github}.

Fix a dimension $d \geq 2$, a desired fractal dimension $\alpha < d$, and
an \textit{innovation parameter} $\theta \geq 0$ which affects the number of clumps
(but does not affect fractal dimension).
Let $B_d$ denote the closed unit ball in $R^d$. 
Build a nested increasing sequence of finite random rooted trees $(\xi_{n}, E_n)_{n \geq 0}$
according to the following rules. 
Take $\xi_0:= \{\mathbf{0} \}$, i.e. the single point at the origin, which serves as the root,
and $E_0:=\emptyset$.

Suppose $n \geq 1$, $\xi_{n-1} \subset B_d$ is a set consisting of $n$ elements,
and $E_{n-1}$ is a collection of directed edges so $(\xi_{n-1},E_{n-1})$ is a tree rooted at $\mathbf{0}$.
To generate $\xi_{n}$, perform two random experiments, independent of each other and of $\xi_{n-1}$:
\begin{enumerate}
\item
Sample $X$ uniformly at random in $B_d$.

\item
Perform a Bernoulli$(\theta/(\theta + n - 1))$ trial (here $0/0 = 1$). 

\end{enumerate}

If the trial is a \textbf{success}, define $\xi_n:=\xi_{n-1} \cup \{X\}$,
declare the parent of $X$ in the tree to be $\mathbf{0}$, and
define $E_n$ to be $E_{n-1}$ together with this extra directed edge. 
In effect we are seeding the process with a new point which need not be close to
a previous one; if $\theta = 0$, this happens only once.

If the trial is a \textbf{failure} (which cannot occur when $n = 1$) compute
\begin{equation} \label{e:mindist}
\Delta:= \min_{\substack{
 x \in \xi_{n-1} \setminus \{\mathbf{0}\} }
} \|X - x\|.
\end{equation}
If $\xi_{n-1}$ were a uniform random sample of $n$ points in $B_d$, this minimum distance $\Delta$
would scale like $n^{-1/d}$. To force points to clump together
we would like the distance of a new point to the closest previous point
to scale like $n^{-1/\alpha}$, where by assumption $1/\alpha > 1/d$.
In the \textbf{smooth minimum distance} formulation, we accept $X$ 
with probability
\begin{equation} \label{e:papp}
e^{-\Delta n^{1/\alpha}}.
\end{equation}
If $X$ is accepted, define $\xi_n:=\xi_{n-1} \cup \{X\}$, and define $E_n$ to be $E_{n-1}$
together with the directed edge $x\to X$, where $x \in \xi_{n-1} \setminus \{\mathbf{0}\}$
 is the minimizing point in (\ref{e:mindist}) (almost surely unique).

If $X$ is not accepted, keep sampling $X$ again in the same way until acceptance occurs.
The conditional law of $\xi_n$ given $\xi_{n-1}$ is
independent of $\xi_0, \xi_1,  \ldots, \xi_{n-2}$, and so $\{\xi_n\}_{n \geq 0}$ is a Markov process.

Except when $n$ is small, the formula (\ref{e:papp}) has the effect of making rejection 
likely when the proposed $(n+1)$-st point is further than $n^{-1/\alpha}$ from any existing
point in $\xi_{n-1} \setminus \{\mathbf{0}\}$. Clumps form, as shown in Figure \ref{f:agorapp}. 
Inserting $\xi_0:= \{\mathbf{0} \}$ serves to create a centering
effect. When $\theta > 0$, the number of distinct seed points, 
like the number of occupied tables in the Chinese restaurant process,
grows according 
\(
 \theta \log{\frac{n+\theta}{\theta}}
\)
after $n$ points, as explained in Pitman \cite{pit}.

\subsection{Hard Threshold Formulation of Agoraphobic Point Process}

Simulation is slow using (\ref{e:papp}), because of
high rejection rates. Here is an inequivalent alternative formulation with similar properties,
whose rejection rate is much lower. 
Replace the smooth acceptance rate(\ref{e:papp}), 
which in (\ref{e:mindist}) is a continuous function
of the sampled point $X$,  by the \textbf{hard threshold}
\[
1_{\{ \Delta \leq n^{-1/\alpha}\} }.
\]
To achieve this efficiently, 
generate a new point which must lie within distance $n^{-1/\alpha}$ of an existing point, as follows:
\begin{enumerate}
\item
Select $z$ uniformly from the $n-1$ points in $\xi_{n-1} \setminus \{\mathbf{0}\}$.
\item
Sample $Y$ uniformly from the ball of radius $n^{-1/\alpha}$, centered at $z$.
\item
Accept $Y$ with probability $1/k$, where
\[
k:= \sum_{x \in \xi_{n-1} \setminus \{\mathbf{0}\} } 1_{\{ \|Y - x\| \leq 1/n\} }.
\]
\item 
If $Y$ is accepted, set $\xi_n:=\xi_{n-1} \cup \{Y\}$, and declare $z$ to be the parent of $Y$ in the random tree.
If $Y$ is not accepted, keep sampling $Y$ again in the same way until acceptance occurs.

\end{enumerate}

\begin{figure}
\caption{\ \textbf{Agoraphobic Point Process:} \textit{
 Two instance of the Hard Threshold Version, 1000 points, with different exponents $1/\alpha$. Hausdorff dimensions of the limiting processes would be
 $10/9$ and $3/2$, respectively.}
} \label{f:agorapp}
\begin{flushleft}
\
\end{flushleft}
\begin{center}
\scalebox{0.35}{\includegraphics{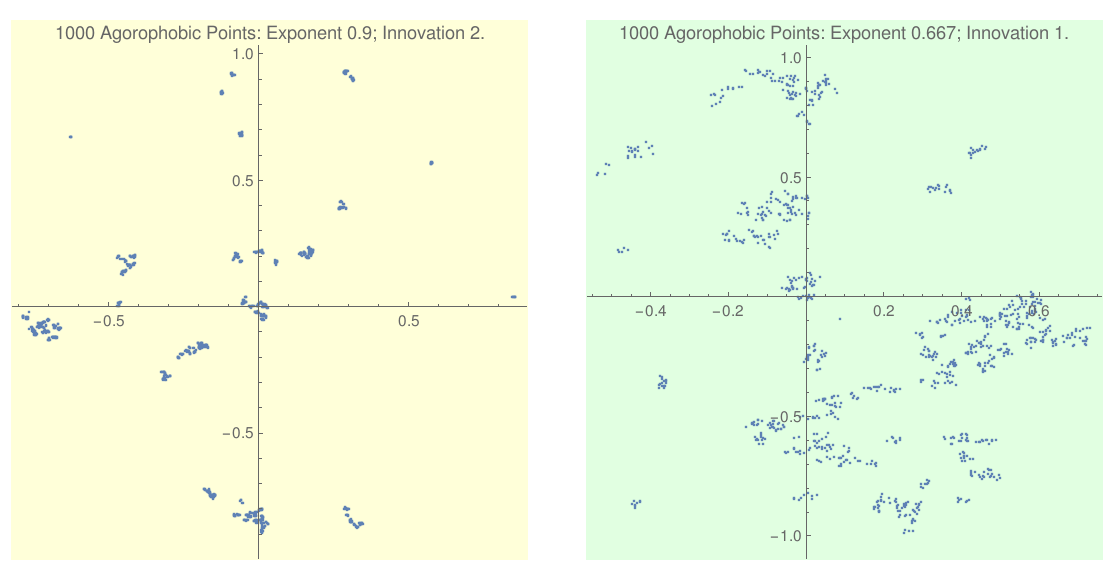}} 
\end{center}
\begin{flushright}
\
\end{flushright}
\end{figure}

Step 3. prevents oversampling in areas which are already dense. 
Not all the points
need lie inside the unit disk. See Figure \ref{f:agorapp}.

\textbf{Acknowledgment: }
The authors thank Dylan Molho (Michigan State Univ.) for analysis of several data sets relevant to this study during the 2017 NSA Summer Program in Operations Research. 
%
%
%
%

\end{document}